\documentclass{article}

\usepackage{amsmath}
\usepackage{amssymb}

\newtheorem{thm}{Theorem}
\newtheorem{conj}{Conjecture}
\newtheorem{lem}{Lemma}

\newtheorem{cor}{Corollary}

\title{Distance between arithmetic progressions and perfect squares}
\author{Tsz Ho Chan}
\date{}

\begin{document}
\maketitle

\begin{abstract}
In this paper, we study how close the terms of a finite arithmetic
progression can get to a perfect square. The answer depends on the
initial term, the common difference and the number of terms in the
arithmetic progression.
\end{abstract}

Many questions in number theory can be phrased as the study of the
``distances" between two sequences of numbers. For instance, we have
the famous conjecture that there are infinitely many primes of the
form $n^2 + 1$. This can be interpreted as saying that the sequence
of prime numbers and the sequence of perfect squares can get within
unit distance from one another infinitely often. Another example is
the conjecture of M. Hall, Jr. stating that
\[
\text{ if } x^2 \neq y^3, \text{ then } |x^2 - y^3| \gg x^{1/2}.
\]
This means that there are significant gaps between the sequence of
perfect squares and the sequence of perfect cubes (apart from the
sequence of sixth powers). In this paper, we are going to consider
the distance between two simplest arithmetic sequences, namely
arithmetic progressions and perfect squares.

More specifically, given integers $a \ge 0$, $d \ge 1$ and $N \ge
1$, we consider the arithmetic progression
\[
\mathcal{A} = \mathcal{A}_{a,d,N} := \{ a, a + d, a + 2d, ..., a + Nd \}.
\]
We are interested in how close the terms of the above arithmetic
progression can get to a perfect square $0^2, 1^2, 2^2, 3^2, ...$.
More precisely, we let
\[
\delta = \delta_{a,d,N} := \mathop{\min_{0 \le n \le N}}_{m \in \mathbb{Z}} |a + nd - m^2|.
\]
To avoid triviality, we require the interval $[a, a + Nd]$ to contain at least one perfect
square. So we need
\[
M := \sqrt{a + Nd} - \sqrt{a} \ge 1 \; \Leftrightarrow \; \sqrt{a + Nd} \ge 1 + \sqrt{a}
\; \Leftrightarrow a \le \Bigl(\frac{Nd - 1}{2}\Bigr)^2
\]
which we will assume throughout the paper. Hence $a \ll N^2 d^2$ in
order for the question to be interesting. Clearly we have the bound
\[
\delta \le \frac{d}{2}.
\]
Let us first give a heuristic investigation. With the above notations, there are about
$M$ perfect squares in the interval $[a, a + Nd]$. Suppose that these squares are
``uniformly distributed" $(\bmod \; d)$. Then we expect that they are spaced about
$\frac{d}{M}$ from one another $(\bmod \; d)$ which implies that one can find a perfect
square within a distance
\[
\frac{d}{M} = \frac{d}{\sqrt{a + Nd} - \sqrt{a}} = \frac{\sqrt{a + Nd} +
\sqrt{a}}{N} \approx \Bigl\{
\begin{array}{ll}
\sqrt{a}/N, & \hbox{ if } a \ge Nd, \\
\sqrt{d}/\sqrt{N}, & \hbox{ if } a \le Nd.
\end{array}
\]
from some $a + n d$ with $0 \le n \le N$. There is a change of
behavior depending on how $a$ compares with $N d$. The above
heuristic makes sense only when $\frac{d}{M} \ge 1$ as the spacing
between integers is at least one.
\begin{equation} \label{Md}
M \le d \Leftrightarrow \sqrt{a + Nd} - \sqrt{a} \le d \Leftrightarrow
N \le d + 2\sqrt{a}
\end{equation}
after some simple algebra. Thus we are led to the following
\begin{conj} \label{conj1}
For any $\epsilon > 0$ and integers $a \ge 0 $, $d \ge 1$, $N \ge 1$ such that $a \le
(\frac{Nd - 1}{2})^2$ and $N \le d + 2\sqrt{a}$,
\[
\delta \ll_\epsilon \frac{\sqrt{a + Nd}}{N} d^\epsilon \hbox{ or equivalently, }
\delta \ll_\epsilon \Bigl\{
\begin{array}{ll}
\frac{\sqrt{a}}{N} d^\epsilon, & \hbox{ if } Nd \le a \le (\frac{Nd - 1}{2})^2, \\
\frac{\sqrt{d}}{\sqrt{N}} d^\epsilon, & \hbox{ if } a \le Nd.
\end{array}
\]
\end{conj}
It may be worthwhile to mention that if $N \le d + 2\sqrt{a}$ and $N
d \le a$, then $N^2 \le Nd + 2N\sqrt{a} \le a + 2N\sqrt{a}$ which
implies $2N^2 \le (\sqrt{a} + N)^2$. Hence $N \ll \sqrt{a}$ and the
ratio $\frac{\sqrt{a}}{N} \gg 1$ in the first half of the
conjecture. Meanwhile if $N \le d + 2\sqrt{a}$ and $a \le Nd$, then
$N \le d + 2\sqrt{Nd}$ which implies $2N \le (\sqrt{d} +
\sqrt{N})^2$. Hence $\sqrt{N} \ll \sqrt{d}$ and the ratio
$\frac{\sqrt{d}}{\sqrt{N}} \gg 1$ in the second half of the
conjecture.

\bigskip

What happens when $N > d + 2\sqrt{a}$? Retracing the steps in
(\ref{Md}), this means that $M > d$. So there are more perfect
squares in the interval $[a, a+Nd]$ than the common difference $d$
of the arithmetic progression. However $m^2 \equiv (m+d)^2 \; (\bmod
\; d)$. So having more than $d$ consecutive perfect squares does not
make them any closer to the arithmetic progression $\mathcal{A}$
than just having $M = d$ perfect squares. Hence by putting in $N = d
+ 2 \sqrt{a}$ in Conjecture \ref{conj1}, we have

\begin{conj} \label{conj2}
For any $\epsilon > 0$ and integers $a \ge 0 $, $d \ge 1$, $N \ge 1$ such that $a \le
(\frac{Nd - 1}{2})^2$ and $N > d + 2\sqrt{a}$,
\[
\delta \ll_\epsilon d^\epsilon.
\]
\end{conj}

Towards Conjectures \ref{conj1} and \ref{conj2}, we have

\begin{thm} \label{thm1}
For $a \le N^2 d / 1800$,
\[
\delta \ll \left\{ \begin{array}{ll}
\frac{a^{1/4} d^{1/2}}{N^{1/2}}, & \hbox{ if } \max(Nd, N^2) \le a \le \frac{N^2 d}
{1800}, \\
 & \\
\frac{d^{3/4}}{N^{1/4}}, & \hbox{ if } N^2 \le a \le Nd, \\
 & \\
d^{1/2}, & \hbox{ if } a \le N^2.
\end{array} \right.
\]
\end{thm}

\begin{thm} \label{thm2}
For $a \le N^{4/3} d^{4/3} / 200$,
\[
\delta \ll \left\{ \begin{array}{ll}
\frac{a^{1/2}}{N^{1/2}}, & \hbox{ if } \max(Nd, N^{2/3} d^{4/3}) \le a \le \frac{N^{4/3}
d^{4/3}}{200}, \\
 & \\
\frac{d}{a^{1/4}}, & \hbox{ if } \max(Nd, d^{4/3}) \le a \le N^{2/3} d^{4/3}.
\end{array} \right.
\]
\end{thm}

\begin{cor} \label{cor1}
For $a \le \frac{N^2 d}{1800}$, $\delta \ll d^{3/4}$.
\end{cor}

Proof of Corollary \ref{cor1}: It follows immediately from Theorem \ref{thm1}. Clearly
the second and third bound in Theorem \ref{thm1} are $\ll d^{3/4}$ while the first bound
$\frac{a^{1/4} d^{1/2}}{N^{1/2}} \ll \frac{N^{1/2} d^{1/4} d^{1/2}}{N^{1/2}} = d^{3/4}$.
Also the three cases cover all the possible ranges for $a \le \frac{N^2 d}{1800}$.

Theorems \ref{thm1} and \ref{thm2} are far from Conjectures
\ref{conj1} and \ref{conj2}. However if we assume a certain
conjectural bound on an average of twisted incomplete Sali\'{e}
sums, we can prove that Conjecture \ref{conj1} is true for a certain
range of $a$, namely
\begin{thm} \label{thm3}
Assume Conjecture \ref{conj3} in section \ref{salie}. If $d$ is odd,
then for all $\epsilon > 0$, there are some constants $C_1, C_2 > 0$
such that
\[
\delta \le C_1 \frac{\sqrt{a}}{N} d^\epsilon \; \hbox{ if } \; N^2
d^{2/3} \le a \le C_2 N^2 d^{1-\epsilon}.
\]
\end{thm}
However, in view of the above Theorems, much is still unknown when
$a$ is big, namely $\max(N^2 d, N^{4/3} d^{4/3}) \ll a \ll N^2 d^2$.

\bigskip

{\bf Some Notations} Throughout the paper, the notations $f(x) =
O(g(x))$, $f(x) \ll g(x)$ and $g(x) \gg f(x)$ are all equivalent to
$|f(x)| \leq C g(x)$ for some constant $C > 0$. Also $f(x) =
O_\lambda(g(x))$, $f(x) \ll_\lambda g(x)$ or $g(x) \gg_\lambda f(x)$
mean that the implicit constant $C$ may depend on $\lambda$.
\section{Proof of Theorems \ref{thm1} and \ref{thm2}}
Our main tool is a result of Huxley [\ref{H}] on integer points close to a curve.
\begin{thm} \label{thm4}
Let $M \ge 12$ be a positive integer, and let $I$ be a closed interval of length $M$
with integer endpoints. Let $f(x)$ be a real function, twice continuously differentiable
on $I$ with
\[
\frac{144}{M^2} \le \frac{\Delta}{C} \le |f''(x)| \le C \Delta \le 1
\]
where $C \ge 1$ and $\Delta \le 1$ are suitable real parameters. Let $\epsilon$ be a
real number with
\[
\max (3 \sqrt{C^3 \Delta}, 6C^2 \sqrt{\frac{2}{M}} ) \le \epsilon \le \frac{1}{4}.
\]
Then there are distinct integers $m_1$, ..., $m_R$ in $I$ with
\[
||f(m_r)|| \le \epsilon \hbox{ for } 1 \le r \le R, \hbox{ and } R \ge \min \Bigl(
\frac{\epsilon^4 M}{2^4 3^4 C^7 \Delta}, \frac{\epsilon M}{2^7 3^2 C^4} \Bigr).
\]
\end{thm}
We also keep using the following simple fact.
\begin{lem} \label{lem1}
For any $x \ge 0$, we can always find a perfect square that is within a distance
$2\sqrt{x}$ from $x$.
\end{lem}

Proof: Observe that $0 \le x - \lfloor \sqrt{x} \rfloor^2 =
(\sqrt{x} - \lfloor \sqrt{x} \rfloor) (\sqrt{x} + \lfloor \sqrt{x}
\rfloor) \le 2 \sqrt{x}$.

\bigskip

Proof of Theorem \ref{thm1}: One can easily check that the theorem
is true when $d = 1$. So we may assume $d \ge 2$ from now on.
Consider $f(x) = \frac{x^2 - a}{d}$. Then $f''(x) = \frac{2}{d} =
\Delta$ and $C = 1$ in the notations of Theorem \ref{thm4}. Without
loss of generality, we may assume that $\sqrt{a+Nd} - \sqrt{a} \ge
14$ for otherwise the theorem is true by picking a large enough
implicit constant. Let $I$ be the closed interval with endpoints
$\lceil \sqrt{a} \rceil$ and $\lfloor \sqrt{a + Nd} \rfloor$. Then
$I$ has length $M = \lfloor \sqrt{a + Nd} \rfloor - \lceil \sqrt{a}
\rceil \ge 12$. For below, we will use
\[
\frac{Nd}{2(\sqrt{a+Nd} + \sqrt{a})} = \frac{\sqrt{a + Nd} - \sqrt{a}}{2} \le M \le
\sqrt{a+Nd} - \sqrt{a} = \frac{Nd}{\sqrt{a+Nd} + \sqrt{a}}.
\]
There are two cases.

\bigskip

Case 1: $a \ge Nd$. Then $\frac{Nd}{5 \sqrt{a}} \le M \le \frac{Nd}{2 \sqrt{a}}$. The
condition $\frac{144}{M^2} \le \frac{\Delta}{C}$ is satisfied when $a \le \frac{N^2
d}{1800}$. There are two subcases.

Subcase 1: $a \le N^2$. Then $\sqrt{\Delta} \gg \frac{1}{\sqrt{M}}$. So we can pick
$\epsilon \ll \sqrt{\Delta} \ll \frac{1}{\sqrt{d}}$ and $m_1 \in I$ with
$||f(m_1)|| \le \epsilon$. This gives
\[
\Big| \frac{m_1^2 - a}{d} - n \Big| \ll \frac{1}{\sqrt{d}} \Rightarrow
|m_1^2 - a - n d| \ll d^{1/2}
\]
for some integer $n$. Since $m_1 \in I$, $0 \le n \le N$. We have the third
bound for Theorem \ref{thm1} when $a \ge Nd$.

Subcase 2: $a \ge N^2$. Then $\sqrt{\Delta} \ll \frac{1}{\sqrt{M}}$. So we can pick
$\epsilon \ll \frac{1}{\sqrt{M}} \ll \frac{a^{1/4}}{N^{1/2} d^{1/2}}$ and $m_1 \in I$ with
$||f(m_1)|| \le \epsilon$. This gives
\[
\Big| \frac{m_1^2 - a}{d} - n \Big| \ll \frac{a^{1/4}}{N^{1/2} d^{1/2}} \Rightarrow
|m_1^2 - a - n d| \ll \frac{a^{1/4} d^{1/2}}{N^{1/2}}
\]
for some integer $0 \le n \le N$, and we have the first bound for Theorem \ref{thm1}.

\bigskip

Case 2: $a \le Nd$. Then $\frac{\sqrt{Nd}}{5} \le M \le \sqrt{Nd}$. The
condition $\frac{144}{M^2} \le \frac{\Delta}{C}$ is satisfied when $N \ge 1800$.
Again there are two subcases.

Subcase 1: $a \le N^2$. Then $\sqrt{\Delta} \gg \frac{1}{\sqrt{M}}$. So we can pick
$\epsilon \ll \sqrt{\Delta} \ll \frac{1}{\sqrt{d}}$ and $m_1 \in I$ with
$||f(m_1)|| \le \epsilon$. This gives
\[
\Big| \frac{m_1^2 - a}{d} - n \Big| \ll \frac{1}{\sqrt{d}} \Rightarrow
|m_1^2 - a - n d| \ll d^{1/2}
\]
for some integer $0 \le n \le N$. We have the third bound for
Theorem \ref{thm1} when $a \le Nd$. Note that if $N < 1800$, by Lemma \ref{lem1},
we can find a perfect square within a distance $2\sqrt{a} \le 2 \sqrt{Nd} \ll d^{1/2}$
from $a$. So we still have the third bound.

Subcase 2: $a \ge N^2$. Then $\sqrt{\Delta} \ll \frac{1}{\sqrt{M}}$. So we can pick
$\epsilon \ll \frac{1}{\sqrt{M}} \ll \frac{a^{1/4}}{N^{1/2} d^{1/2}}$ and $m_1 \in I$ with
$||f(m_1)|| \le \epsilon$. This gives
\[
\Big| \frac{m_1^2 - a}{d} - n \Big| \ll \frac{1}{N^{1/4} d^{1/4}} \Rightarrow
|m_1^2 - a - n d| \ll \frac{d^{3/4}}{N^{1/4}}
\]
for some integer $0 \le n \le N$, and we have the second bound for Theorem \ref{thm1}.
Note that if $N < 1800$, by Lemma \ref{lem1}, we can find a perfect square within a
distance $2\sqrt{a} \le 2 \sqrt{Nd} \ll d^{1/2} \ll \frac{d^{3/4}}{N^{1/4}}$ from $a$.
So we still have the second bound.

\bigskip

Proof of Theorem \ref{thm2}: We apply Theorem \ref{thm4} to the
inverse function $g(x) = f^{-1}(x) = \sqrt{a + x d}$ and $I$ is the
interval $[0,N]$ with length $M = N$. Then
\[
g''(x) = - \frac{d^2}{4 (a + x d)^{3/2}}.
\]
We shall focus on the case $a \ge Nd$ only as the other case is covered by Theorem
\ref{thm1}. We have
\[
\frac{d^2}{16 a^{3/2}} \le |g''(x)| \le \frac{d^2}{4 a^{3/2}}.
\]
So $\Delta = \frac{d^2}{8 a^{3/2}}$ and $C = 2$ in the notations of Theorem
\ref{thm1}. The condition $C \Delta \le 1$ is equivalent to $\frac{d^{4/3}}{2^{4/3}}
\le a$. Meanwhile the condition $\frac{144}{M^2} \le \frac{\Delta}{C}$ is satisfied
when $a \le \frac{N^{4/3} d^{4/3}}{200}$. There are two subcases.

Subcase 1: $a \le N^{2/3} d^{4/3}$. Then $\sqrt{\Delta} \gg \frac{1}{\sqrt{M}}$.
So we can pick $\epsilon \ll \sqrt{\Delta} \ll \frac{d}{a^{3/4}}$ and $n \in I$ with
$||g(n)|| \le \epsilon$. This gives, for some integer $m$,
\[
|\sqrt{a + n d} - m| \ll \frac{d}{a^{3/4}} \Rightarrow |a + nd - m^2| =
|\sqrt{a + n d} - m| |\sqrt{a + n d} + m| \ll \frac{d}{a^{1/4}}
\]
as $a \gg d^{4/3}$ implies $\frac{d}{a^{3/4}} \ll a^{1/2}$. This gives the second
bound for Theorem \ref{thm2}.

Subcase 2: $a \ge N^{2/3} d^{4/3}$. Then $\sqrt{\Delta} \ll \frac{1}{\sqrt{M}}$.
So we can pick $\epsilon \ll \frac{1}{\sqrt{M}} = \frac{1}{\sqrt{N}}$ and $n \in I$ with
$||g(n)|| \le \epsilon$. This gives, for some integer $m$,
\[
|\sqrt{a + n d} - m| \ll \frac{1}{N^{1/2}} \Rightarrow |a + nd - m^2| =
|\sqrt{a + n d} - m| |\sqrt{a + n d} + m| \ll \frac{a^{1/2}}{N^{1/2}}.
\]
This gives the first bound for Theorem \ref{thm2}.
\section{Twisted Sali\'{e} Sum and Theorem \ref{thm3}} \label{salie}

The same technique in [\ref{FI}] and [\ref{C}] for the studies of
$n^2 \alpha \pmod 1$ and short intervals containing almost squares
or sums of two squares can be used here. We recall the following
conjectural bound on a certain average of twisted incomplete
Sali\'{e} sums.
\begin{conj} \label{conj3}
Let $a$, $q$ be integers with $q \ge 2$ and $(a,q) = 1$ and $q$ is not a perfect square.
Let $H, K \ge 1$ and $\lambda$, $\mu$ be any real numbers. Then, for any
$\epsilon > 0$,
\begin{align*}
\mathop{\sum_{1 \le h \le H}}_{(h,q) = 1} e(\lambda h) \sum_{0 \le k < K}& e(\mu k)
\Bigl(\frac{h}{q} \Bigr) e \Bigl(\frac{a \overline{h} k^2}{q} \Bigr) \\
\ll_\epsilon& (H^{1/2} K^{1/2} + H^{3/4} + K + q^{-1/2} HK + q^{-1/2} K^2) q^\epsilon.
\end{align*}
\end{conj}
A consequence of the above conjecture is the following
\begin{lem} \label{lem2}
Assume Conjecture \ref{conj3}. Let $q$ be an odd number, $H, K \ge 1$ and $\lambda$,
$\mu$ be any real numbers. We have
\begin{align*}
\sum_{0 \le k < K} e(\mu k) & \mathop{\sum_{1 \le h \le H}}_{(h,q) = 1} e(\lambda h)
G(h,\pm k; q) \\
\ll_\epsilon & (q^{1/2} H^{1/2} K^{1/2} + q^{1/2} H^{3/4} + q^{1/2} K + HK + K^2)
q^\epsilon
\end{align*}
where
\[
G(a,b;q) := \sum_{n \; (\bmod \; q)} e \Bigl( \frac{a n^2 + b n}{q} \Bigr)
\]
is the Gauss sum.
\end{lem}

Proof: It is Lemma 5.1 in [\ref{C}].

\bigskip

Proof of Theorem \ref{thm3}: It is very similar to the proof of Theorem 1.4 in [\ref{C}].
Suppose that $N$ and $d$ are sufficiently large and $a \ge \max(Nd, N^2)$. Let
$A := \lfloor \sqrt{a + \frac{N}{2} d} \rfloor$. We restrict
our attention to the arithmetic progression $a + n d$ where $\frac{N}{4} \le n \le
\frac{3N}{4}$. Suppose $a + n d = (A + m)^2$. Then
\[
m = \sqrt{a + n d} - \sqrt{a + \frac{N}{2} d} + O(1) = \frac{(n - N/2) d}{\sqrt{a + n d}
+ \sqrt{a + N d / 2}} + O(1).
\]
Let $M = \lfloor \frac{Nd}{12 \sqrt{a}} \rfloor$. If $m \in [-M, M]$, then $n \in
[\frac{N}{4}, \frac{3N}{4}]$. Our goal is to find some $m \in [-M, M]$ such that
\[
\Big\| \frac{(A + m)^2 - a}{d} \Big\| < \Delta \Rightarrow \Big| \frac{(A + m)^2 - a}{d}
- n \Big| < \Delta \Rightarrow |(A + m)^2 - (a + n d)| < \Delta d
\]
for some integer $n \in [\frac{N}{4}, \frac{3N}{4}]$. For any $0 < \Delta < 1/2$, define
\[
S := \# \Bigl\{m : |m| \le M, \; \Big\| \frac{(A + m)^2 - a}{d} \Big\| < \Delta \Bigr\}.
\]
Let
\[
f(x) := \left\{ \begin{array}{ll}
1 - |x|, & \hbox{ if } -1 \le x \le 1, \\
0, & \hbox{ otherwise,}
\end{array} \right.
\; \;
t(x) := \left\{ \begin{array}{ll}
1 - |x/ \Delta|, & \hbox{ if } -\Delta \le x \le \Delta, \\
0, & \hbox{ otherwise,}
\end{array} \right.
\]
\[
g(x) := \sum_{n = -\infty}^{\infty} t(x - n), \; \; g_{\lambda}(x) = g(x - \lambda).
\]
Thus $f(x)$ can be used to detect numbers between $-1$ and $1$, and $g_\lambda(x)$ can be
used to detect numbers whose fractional part is between $\lambda - \Delta$ and $\lambda +
\Delta$. In particular with $\lambda = a/d$,
\[
S \ge \sum_{m = -\infty}^{\infty} f \Bigl(\frac{m}{M} \Bigr) g_\lambda \Bigl(
\frac{(A+m)^2}{d} \Bigr).
\]
Now $g_\lambda$ has Fourier expansion
\[
g_\lambda(x) = \sum_{h = -\infty}^{\infty} c(h) e(-\lambda h) e(h x)
\]
where
\[
e(u) = e^{2 \pi i u}, \; \; \hbox{ and } \; \; c(h) = \Delta \Bigl( \frac{\sin
\pi \Delta h}{\pi \Delta h} \Bigr)^2 \hbox{ when } h \neq 0, \; \; c(0) = \Delta.
\]
Therefore
\[
S = \Delta M + \sum_{h \neq 0} c(h) e(-\lambda h) \sum_{m = -\infty}^{\infty}
e \Bigl( \frac{h (A+m)^2}{d} \Bigr) f \Bigl( \frac{m}{M} \Bigr) =: \Delta M + R.
\]
By Poisson summation, the inner sum over $m$ in $R$
\begin{align*}
=& \sum_{t \; (\bmod \; d)} e \Bigl( \frac{h (A+t)^2}{d} \Bigr) \sum_{m \equiv t \;
(\bmod \; d)} f \Bigl( \frac{m}{M} \Bigr) \\
=& \sum_{t \; (\bmod \; d)} e \Bigl( \frac{h (A+t)^2}{d} \Bigr) \frac{M}{d} \sum_{k =
-\infty}^{\infty} \hat{f} \Bigl(\frac{Mk}{d}\Bigr) e \Bigl(\frac{t k}{d} \Bigr) \\
=& \frac{M}{d} \sum_{k = -\infty}^{\infty} \hat{f} \Bigl(\frac{Mk}{d}\Bigr)
\sum_{t \; (\bmod \; d)} e \Bigl(\frac{h(A+t)^2 + k t}{d} \Bigr) \\
=& \frac{M}{d} \sum_{k = -\infty}^{\infty} \hat{f} \Bigl(\frac{Mk}{d}\Bigr)
e \Bigl(- \frac{k A}{d} \Bigr) \sum_{t' \; (\bmod \; d)} e \Bigl(\frac{h t'^2 + k t'}{d}
\Bigr) = \frac{M}{d} \sum_{k = -\infty}^{\infty} \hat{f} \Bigl(\frac{Mk}{d}\Bigr)
e \Bigl(- \frac{k A}{d} \Bigr) G(h, k; d)
\end{align*}
where $\hat{f}(y) = (\frac{\sin \pi y}{\pi y})^2$ is the Fourier transform of $f$ and
$G(h,k;d)$ is the Gauss sum. Thus
\begin{align*}
R =& \frac{M}{d} \sum_{k = -\infty}^{\infty} \Bigl( \frac{\sin(\pi k M /d)}{\pi k M / d}
\Bigr)^2 e \Bigl(- \frac{k A}{d} \Bigr) \sum_{h \neq 0} c(h) e(-\lambda h) G(h,k;d) \\
=& \frac{M}{d} \Bigl( \sum_{|k| < d/M} + \sum_{d/M \le |k| < L} + \sum_{|k| \ge L} \Bigr)
\end{align*}
where $L = (d/M)^\sigma$ for some parameter $\sigma > 1$ to be chosen later. At this point,
the proof proceeds almost identically as that of Lemma 6.1 in [\ref{C}] by estimating
$\sum_{|k| \ge L}$ trivially and the other two sums using Lemma \ref{lem2}. Apart from
different choices of letters for the variables, the only difference in the argument is when applying partial summation, we have
\[
\sum_{0 \le k < d/M} = \int_{0}^{d/M} \Bigl( \frac{\sin(\pi u M
/d)}{\pi u M / d} \Bigr)^2 d \Bigl( \sum_{0 \le k < u} e \Bigl(-
\frac{k A}{d} \Bigr) \sum_{h \neq 0} c(h) e(-\lambda h) G(h,k;d)
\Bigr)
\]
for example (in contrast with equation (13) in [\ref{C}]). Also we shall keep $\sigma$
(i.e. $L$) in our bound instead of choosing $\sigma = 5$ in [\ref{C}]. Consequently, we have
\[
R \ll_\epsilon d^\epsilon \log \frac{d}{\Delta} \log L \Bigl[M^{1/2} \Delta^{1/2} +
\frac{M \Delta^{1/4}}{d^{1/2}} + \Delta d^{1/2} + 1 + \frac{d \Delta}{M} +
\frac{d^2}{\Delta M L} \Bigr].
\]
Thus for $\Delta \ge \frac{1}{d}$ and $a \le N^2 d^2$, there exists some constant
$C_\epsilon > 0$ such that
\[
R \le C_\epsilon \sigma d^{2 \epsilon} \Bigl[M^{1/2} \Delta^{1/2} +
\frac{M \Delta^{1/4}}{d^{1/2}} + \Delta d^{1/2} + 1 + \frac{d \Delta}{M} +
\frac{d^2}{\Delta M L} \Bigr].
\]
We are going to make each of the six pieces of the above bound for $R$ less than $\Delta M
/ 12$. Then it follows that $S > \Delta M / 2 > 0$. The six pieces
of requirements are achieved if (i) $\Delta \gg \frac{\sqrt{a}}{N d}
d^{4 \epsilon}$, (ii) $\Delta \gg \frac{d^{3\epsilon}}{d^{2/3}}$, (iii) $a \ll N^2 d^{1 -
4\epsilon}$, (iv) $\Delta \gg \frac{\sqrt{a}}{Nd} d^{2\epsilon}$, (v) $a \ll N^2 d^{1 -
2\epsilon}$, and (vi) $\Delta \gg \frac{N^{\sigma/2 - 1}}{a^{\sigma/4 - 1/2}} d^\epsilon$
respectively where the implicit constants may depend on $\epsilon$ and $\sigma$. Since
$\epsilon$ is arbitrary, the above can be summarized as
\[
\Delta \gg \frac{\sqrt{a}}{N d} d^{\epsilon}, \; \; \Delta \gg \frac{1}{d^{2/3 -
\epsilon}}, \; \; a \ll N^2 d^{1-\epsilon}, \; \; \Delta \gg \frac{N^{\sigma/2 - 1}}
{a^{\sigma/4 - 1/2}} d^\epsilon.
\]
Now consider $N^2 d^{2/3} \le a \ll_\epsilon N^2 d^{1-\epsilon}$. The first lower bound
on $\Delta$ is the biggest by taking $\sigma$ large enough. Therefore there exist some
$m \in [-M,M]$ and $n \in [\frac{N}{4}, \frac{3N}{4}]$ such that
\[
|(A+m)^2 - (a + n d)| < \Delta d \ll_\epsilon \frac{\sqrt{a}}{N} d^{\epsilon}
\]
which proves Theorem \ref{thm3}.

Tsz Ho Chan \\
Department of Mathematical Sciences \\
University of Memphis \\
Memphis, TN 38152 \\
U.S.A. \\
tchan@memphis.edu \\

\end{document}


\maketitle


\section{Introduction and Main Result}

Before, nothing was known in the range $\max(N^2 d, N^{4/3} d^{4/3}) \ll a \ll N^2 d^2$. In this note, we are going to answer this, namely,
\begin{thm} \label{thm3}
For infinitely many positive integers $N$, $d \ge 30$ and $25 N^2 d \le a \le N^2 d^2$,
\[
\delta \ge \frac{d}{1800}.
\]
\end{thm}
Thus, the minimum distance has the same order of magnitude as the trivial bound $d/2$ in this range of $a$. It is also interesting to note the sudden change in behavior from $a \le N^2 d / 1800$ to $a \ge 25 N^2 d$. So the next challenge would be to determine the exact moment and nature of transition from order $d^{3/4}$ to $d$ when $N^2 d \ll a \ll N^2 d$.

%
%
\section{Proof of Theorem \ref{thm3}}

Suppose $a$ is not a perfect square. We want to find $0 \le n \le N$ so that $\sqrt{a + n d}$ is close to an integer. By Taylor's expansion with remainder,
\begin{equation} \label{taylor}
\sqrt{a + n d} = \sqrt{a} \sqrt{1 + \frac{n d}{a}} = \sqrt{a} \Bigl[1 + \frac{n d}{2 a} + \frac{\epsilon_0 n^2 d^2}{a^2} \Bigr] = \sqrt{a} + \frac{n d}{2 \sqrt{a}} + \frac{\epsilon_0 n^2 d^2}{a^{3/2}}
\end{equation}
for some $|\epsilon_0| \le 1/8$. Denote
\[
A := \lfloor \sqrt{a} \rfloor.
\]
Since $A \le \sqrt{a} + \frac{n d}{2 \sqrt{a}} \le A + \frac{N d}{2 \sqrt{a}} + 1$, for each integer $0 \le k \le \frac{N d}{2 \sqrt{a}} + 1$, there is at most one integer $n_k \le N$ such that
\[
\Big| \sqrt{a} + \frac{n_k d}{2 \sqrt{a}} - (A + k) \Big| \le \frac{d}{4 \sqrt{a}}.
\]
For all integers $n$ other than these $n_i$'s, we must have
\begin{equation} \label{other}
\Big| \sqrt{a} + \frac{n d}{2 \sqrt{a}} - m \Big| > \frac{d}{4 \sqrt{a}}
\end{equation}
for any integer $m$. Now, we restrict our attention to specific $a$'s and $d$'s so that the fractional part of $2 \sqrt{a} / d$ is close to $1/3$. Consider
\[
d = 2(9 d' + 1) \text{ and } a = \frac{(9 d' + 1)^2 (9 X + 1)^2 - 1}{9}
\]
for some integers $d'$ and $X$. With these special forms of $a$ and $d$, we have
\begin{equation} \label{e1}
\frac{2 \sqrt{a}}{d} = \Big(3 X + \frac{1}{3} \Bigr) \sqrt{1 - \frac{1}{(9 d' + 1)^2 (9 X + 1)^2}} = 3 X + \frac{1}{3} + \frac{\epsilon_1}{d \sqrt{a}}
\end{equation}
for some $-3 < \epsilon_1 < 0$. Also,
\begin{equation} \label{e2}
\sqrt{a} = (9 d' + 1) \Bigl(3 X + \frac{1}{3} \Bigr) \sqrt{1 - \frac{1}{(9 d' + 1)^2 (9 X + 1)^2}} = (27 d' X + 3 X + 3 d') + \frac{1}{3} + \frac{\epsilon_2}{\sqrt{a}}
\end{equation}
for some $-1 < \epsilon_2 < 0$. In particular, $A = 27 d' X + 3 X + 3 d'$. Combining (\ref{e1}) and (\ref{e2}), we have
\begin{equation} \label{e3}
\frac{2 \sqrt{a}}{d} (\sqrt{a} - A) = X + \frac{1}{9} + \frac{\epsilon_3}{d}
\end{equation}
for some $-3 < \epsilon_3 < 0$. Therefore, by (\ref{e1}) and (\ref{e3}), the fractional part of
\[
\frac{2 \sqrt{a}}{d} k - \frac{2 \sqrt{a}}{d} (\sqrt{a} - A)
\]
is at least $\frac{1}{9} + \frac{\epsilon_4}{d}$ with $-3 < \epsilon_4 < 3$ for any integer $k \le \frac{N d}{2 \sqrt{a}} + 1$ since $N^2 d \le a \le N^2 d^2$. Consequently,
\[
\Big| \sqrt{a} + \frac{n_k d}{2 \sqrt{a}} - (A + k) \Big| \ge \Bigl( \frac{1}{9} + \frac{\epsilon_4}{d} \Bigr) \frac{d}{2 \sqrt{a}}.
\]
Combining this with (\ref{other}), we have
\[
\Big| \sqrt{a} + \frac{n d}{2 \sqrt{a}} - m \Big| > \frac{d}{180 \sqrt{a}}
\]
for any non-negative integers $m$ and $n$ as long as $d = 2(9d'+1) \ge 30$. This together with (\ref{taylor}) gives
\[
|\sqrt{a + n d} - m| > \frac{d}{1800 \sqrt{a}}
\]
for any non-negative integers $m$ and $n$ if $a \ge 25 N^2 d$. This, in turn, implies that
\[
|(a + n d) - m^2| > \frac{d}{1800}
\]
and hence our theorem.

%
%